\documentclass[preprint]{elsarticle}
\biboptions{sort&compress}

\newtheorem{thm}{Theorem}

\newdefinition{rmk}[thm]{Remark}
\newdefinition{crl}[thm]{Corollary}
\newproof{pf}{Proof}
\newproof{pot}{Proof of Theorem \ref{thm2}}
\newdefinition{definition}[thm]{Definition}

\usepackage{amssymb}
\usepackage{amsmath}	

\begin{document}

\begin{frontmatter}

\title{Weakly nonlinear hyperbolic differential equation in Hilbert space}

\author{Oleksandr Pokutnyi \footnote{The author acknowledges the financial support from the grant ``Boundary-Value Problems and Impulse Perturbations of Nonlinear Evolution Equations in Infinite-Dimensional Spaces" $3M2022$ (Reg. No. 0122U002463).}}

\ead{lenasas@gmail.com, alex_poker@imath.kiev.ua}

\begin{abstract}
We consider nonlinear perturbations of the hyperbolic equation in the Hilbert space. Necessary and sufficient conditions for the existence of solutions of boundary-value problem for the corresponding equation and iterative procedures for their finding are obtained in the case when the operator in linear part of the problem hasn't inverse and can have nonclosed set of values. As an application we consider boundary-value problems for countable systems of differential equations and van der Pol equation in a separable Hilbert space.
\end{abstract}


\end{frontmatter}

\section{Introduction}
It is well known that the operator approach is a powerful tool for the investigation of boundary-value problems. We can write general boundary-value problems in the following operator form
\begin{equation} \label{war:1}
Ly(t, \varepsilon) = f(t) + \varepsilon Z(y(t, \varepsilon), t, \varepsilon),
\end{equation}
\begin{equation} \label{war:2}
l y(\cdot, \varepsilon) = \alpha + \varepsilon J(y(\cdot, \varepsilon), \cdot, \varepsilon)
\end{equation}
where $L$ and $l$ are linear and bounded operators in some function spaces, $Z, J$ are nonlinear operator-function and vector-functional. A numerous number of articles are devoted to the case when the operator $L$ is a differential operator, operator with delay, with impulse actions with values in finite-dimensional space etc. As a rule (in the finite dimensional case) we have the following cases:
boundary-value problem (\ref{war:1}), (\ref{war:2}) is (regular) everywhere solvable  (according to the classification of Krein S.G. \cite{Krein-1967} (for any right hand side $f(t), \alpha$ there is a unique solution of (\ref{war:1}), (\ref{war:2})), Fredholm with zero index or Noetherian (Fredholm with nonzero index).  Using the theory of generalized inverse operators, it became possible to study boundary value problems in the irregular case when the corresponding problem is not solvable for any right-hand sides and may not have a unique solution \cite{Boi-Sam-2016}. In the given paper we consider infinite-dimensional case when the hyperbolic operator in the linear part can be or normally solvable (with the closed set of values) or can have nonclosed set of values \cite{BoiPok} (strong generalized normally-solvable). As it is known hyperbolic systems of equations are a class of wave equations that can be used to describe various processes in nature \cite{Evans}, \cite{Gabov}. We consider the nonlinear boundary-value problems for the class of abstract hyperbolic equations in the Banach and Hilbert spaces. The theory of differential equations in such spaces occupies an important place in modern mathematics. A great number of fundamental results have been obtained in this area \cite{Krein-1967,Gorbachuk-1984,Fattor-1985,Showalt-1994,Shklyar-1997,Prato-2002,Gav-Mak-Vas-2010,LyashNomir-2012,Diagana-2018}.
Nevertheless, research on various aspects of the qualitative theory of such equations remains relevant \cite{Swiech-1994,Artamon-2003,BaiLiGe-2005,Modanl-Akgul-2017,Gil-2018,Kumar-Mus-Sak-2018,Gav-Mak-May-2021}. Such questions include, in particular, the study of the existence and construction of solutions of boundary-value problems for differential equations with unbounded operator coefficients in Banach and Hilbert spaces. In \cite{Boi-Pok-2019}, a weakly perturbed linear boundary-value problem for a second-order differential equation of hyperbolic type in a separable Hilbert space is considered. The criterion for the existence of solutions of the generating unperturbed problem is obtained and sufficient conditions for the appearance of solutions of the weakly perturbed problem are established under condition that the generating problem has no solutions. In the presented paper, we consider one of the possible approaches for finding necessary and sufficient conditions for the existence of solutions to weakly nonlinear second-order differential equation of hyperbolic type in the Hilbert space and propose an algorithm for finding these solutions. Moreover, the  suggested methodology   also works in the case when the corresponding operator in the linear part has a nonclosed set of values. As an application we consider abstract evolution van der Pol equation in the separable Hilbert space. It is well known that van der Pol equation \cite{Pol-Mark-1920,Pol-Mark-1926} used in modeling the work of the heart \cite{Pol-Mark-1928}, demultiplication of frequency \cite{Van_Der},
in the study of hysteresis \cite{App}, chaos and sincronisation \cite{Mina}, \cite{Zheng}, neural interactions \cite{Wilson-Cowan-1972,Kawahara-1980}, oscillatory movements of human limbs \cite{Beek-Sch-Mor-Sim-Tur-1995}, plate interactions in a geological fault \cite{Cartw-Egu-Her-Gar-Piro-1999}, rotations of the automated system \cite{Veskos-Demiris-2005}, vocal cord modeling \cite{Lucero-Schoentgen-2013}, in the theory of nonlinear oscillations \cite{Levsh}. Problems for the van der Pol equation find its application in singularity theory \cite{Lev}. In \cite{Huang} variational methods and He polynomials are used for investigation of the corresponding problem. Van der Pol equation is also used in algebraic geometry \cite{Odani}, fractional calculus \cite{Mina}, \cite{Mina1} and biological applications \cite{Kapl-2008}.  In paper\cite{Chunru}, bifurcations for a system of three oscillators connected to the van der Pol equation are considered. Presented in the article abstract van der Pol equation model can describe the set of coupled oscillators (in general, their number can be countable).

\section{Statement of the problem}
We consider a boundary-value problem for the abstract hyperbolic equation in the separable Hilbert space $\mathcal{H}$:
\begin{align}
&y^{\prime\prime}(t,\varepsilon) + A(t)y(t,\varepsilon) = \varepsilon Z(t,y(t,\varepsilon), y^{\prime}(t,\varepsilon)) + f(t), \label{Pok-hyperbolic-eq:123} \\
&l(y(\cdot, \varepsilon), y^{\prime}(\cdot, \varepsilon)) = \alpha, \label{Pok-hyperbolic-eq:124}
\end{align}
where $y \in C^2(J,\mathcal{H})$, $J=[0,w]\subset\mathbb{R}$, the closed strongly continuous operator-valued function $A(t)$ acts from $J$  into the dense domain $D=D(A(t))\subset\mathcal{H}$, which is independent from $t$, $Z(t,y(t,\varepsilon), y^{\prime}(t,\varepsilon))$ is continuous with respect to first component and nonlinear,  strongly differentiable according to Fr\'{e}chet in a certain neighborhood of the generating solution with respect to second and third components of the operator-valued function
\begin{equation*}
Z(\cdot,y(t,\varepsilon), y^{\prime}(t,\varepsilon)) \in C(J,\mathcal{H}),
\end{equation*}
\begin{equation*}
Z(t,\cdot,y^{\prime}(t,\varepsilon)) \in C^1[\|y-y_0\|\leq q], \quad
Z(t,y(t,\varepsilon),\cdot) \in C^1[\|y-y_0\|\leq q],
\end{equation*}
where $q$ is sufficiently small constant, $\varepsilon<<1$ is a small parameter, vector-function $f \in C(J,\mathcal{H})$, $l:C^2[J,\mathcal{H}] \times C^{1}[J, \mathcal{H}] \rightarrow\mathcal{H}_{1}$ is bounded linear vector functional, $\alpha \in \mathcal{H}_{1}$. Thus, we consider an abstract hyperbolic equation according to the classification of S.~G.~Krein's \cite{Krein-1967}.

We investigate the questions of finding conditions for the existence and effective construction of boundary-value problem (\ref{Pok-hyperbolic-eq:123}), (\ref{Pok-hyperbolic-eq:124}) solutions, which for $\varepsilon = 0$ turns into the solution $y_0(t)$ of the generating boundary-value problem
\begin{equation} \label{Pok-hyperbolic-eq:3}
	y^{\prime\prime}_{0}(t) + A(t)y_{0}(t) = f(t), \quad l(y_{0}(\cdot), y^{\prime}_{0}(\cdot)) = \alpha.
\end{equation}
These solutions $y_0(t)$ will be called generating solutions of the boundary value problem (\ref{Pok-hyperbolic-eq:123}), (\ref{Pok-hyperbolic-eq:124}).

\section{Criterion of solvability of the generating boundary-value problem}  We give the criterion for the existence of solutions of the boundary value problem (\ref{Pok-hyperbolic-eq:3}), obtained in \cite{Boi-Pok-2019}.

To do this, we make the following change in variables in the boundary-value problem (\ref{Pok-hyperbolic-eq:3})
\begin{equation*}
x_{1}^{0}(t) = y_{0}(t), \quad x_{2}^{0}(t) = y_{0}^{\prime}(t), \quad x_{0}(t) = {\rm col} \begin{pmatrix}
	x_{1}^{0}(t), & x_{2}^{0}(t)
\end{pmatrix}
\end{equation*}
and rewrite it as a boundary-value problem for the operator system
\begin{align}
& x_{0}^{\prime}(t) = B(t)x_{0}(t) + g(t), \label{Pok-hyperbolic-eq:4} \\
&lx_{0}(\cdot) = \alpha, \label{Pok-hyperbolic-eq:04}
\end{align}
where
\begin{equation}
 B(t) = \begin{pmatrix}
   O & I  \\
   -A(t) & O
 \end{pmatrix}, \quad g(t) = {\rm col} \begin{pmatrix}
	0, & f(t)
\end{pmatrix},
\end{equation}
$O$, $I$ are the zero operator and the identity operator, respectively, in space $\mathcal{H}$.
Note that other substitutions of variables are also possible (see, e.g., \cite{Rouche-1980}).
Denote by $U(t)$ the evolution operator of a homogeneous system
\begin{equation*}
U^{\prime}(t) = B(t)U(t), \quad U(0) = I.
\end{equation*}
Then the set of solutions of (\ref{Pok-hyperbolic-eq:4}) has the form
\begin{equation}\label{Pok-hyperbolic-eq:004}
x_{0}(t, c) = U(t)c + \int\limits_{0}^{t}U(t)U^{-1}(\tau)g(\tau)d\tau.
\end{equation}
Substituting (\ref{Pok-hyperbolic-eq:004}) in the boundary condition (\ref{Pok-hyperbolic-eq:04}), we obtain the following
operator equation
\begin{equation} \label{Pok-hyperbolic-eq:5}
Qc = g_{1},
\end{equation}
where
\begin{equation*}
 Q = lU(\cdot) : \mathcal{H} \rightarrow \mathcal{H}_{1}, \quad g_{1} = \alpha - l\int\limits_{0}^{\cdot}U(\cdot)U^{-1}(\tau)g(\tau)d\tau.
\end{equation*}

It is known (see \cite{Boi-Pok-2019}), under different conditions on the right-hand side of $g_1$, the equation (\ref{Pok-hyperbolic-eq:5}) has three types of solutions: 1) classical generalized solutions; 2) strong generalized solutions; 3) strong pseudosolutions. Here is some theoretical information about these solutions that we will need in the future.

First, we consider the case when the set of values of the operator $Q$ is closed: ($R(Q)=\overline{R(Q)}$).
Then operator equation (\ref{Pok-hyperbolic-eq:5}) is solvable if and only if the element $g_1$ satisfy condition
$g_{1} \in R(Q)$ or $\mathcal{P}_{N(Q^{*})}g_{1} = 0$ or $\mathcal{P}_{\mathcal{H}_{Q}}g_1 = 0$ ($\mathcal{H} = \mathcal{H}_{Q} \oplus R(Q) = (\mathcal{H} \ominus R(Q)) \oplus R(Q)$)\cite{Boi-Sam-2016}.
Here $\mathcal{P}_{N(Q^{*})}$  is the orthoprojector onto the cokernel of the operator $Q$ and $\mathcal{P}_{R(Q)}$ is the orthoprojector onto the subspace $\mathcal{H}_{Q}$.
The set of solutions of equation (\ref{Pok-hyperbolic-eq:5}) has the form
\begin{equation*}
c = Q^{+}g_{1} + \mathcal{P}_{N(Q)}\overline{c} \quad \forall \overline{c} \in \mathcal{H},
\end{equation*}
where $\mathcal{P}_{N(Q)}$ is the orthoprojector onto the kernel of the operator $Q$. Such solutions are called classical generalized solutions.

Let us now we consider the case when the set of values of the operator $Q$ is not closed ($R(Q)\neq\overline{R(Q)}$) and $g_{1} \in \overline{R(Q)}$. In this case, the operator $Q$ can be extended to the operator $\overline{Q}$ with a closed set of values ($R\left(\overline{Q}\right)=\overline{R\left(\overline{Q}\right)}$). Indeed, since the operator $Q$ is linear and bounded, spaces $\mathcal{H}$ and $\mathcal{H}_{1}$ can be represented in the form of a direct sums of subspaces:
\begin{equation*}
\mathcal{H} = N(Q) \oplus X, \quad \mathcal{H}_{1} = \overline{R(Q)}
\oplus Y,
\end{equation*}
where  $X = N(Q)^{\bot}$,  $Y = \overline{R(Q)}^{\bot}$. Let
$\mathcal{H}_{2} = \mathcal{H}/N(Q)$ is quotient space of $\mathcal{H}$ and let $\mathcal{P}_{\overline{R(Q)}}$ and $\mathcal{P}_{N(Q)}$ are orthoprojectors onto $\overline{R(Q)}$ and $N(Q)$, respectively. Then the operator
\begin{equation*}
\mathcal{Q} = \mathcal{P}_{\overline{R(Q)}}Qj^{-1}p : X
\rightarrow R(Q) \subset \overline{R(Q)}
\end{equation*}
is linear, continuous, and injective. Here
\begin{equation*}
p : X \rightarrow \mathcal{H}_{2},  \quad j : \mathcal{H}
\rightarrow \mathcal{H}_{2}
\end{equation*}
are continuous bijection and projection, respectively. The triple $(\mathcal{H}, \mathcal{H}_{2}, j)$ is a locally trivial bundle with a typical fiber
$ \mathcal{P}_{N(Q)}\mathcal{H}$. In this case, we can define a strong generalized solution of the equation
\begin{equation} \label{Pok-hyperbolic-eq:1}
Qc = g_{1}, \quad c \in X.
\end{equation}
We complete the space $X$ by the norm $||x||_{\overline{X}} = ||\mathcal{Q}x||_{F}$ \cite{LyashNomir-2012}, where $F =
\overline{R(Q)}$. Then the extended operator $\overline{\mathcal{Q}}$
\begin{equation*}
\overline{\mathcal{Q}} :  \overline{X}
\rightarrow R(Q) \subset \overline{R(Q)}
\end{equation*}
is a homeomorphism of $\overline{X}$ and $\overline{R(Q)}$. By the construction of a strong generalized
solution, the equation
\begin{equation*}
\overline{\mathcal{Q}}\overline{c} = g_{1}
\end{equation*}
has a unique solution $\overline{\mathcal{Q}}^{-1}g_{1}$, which is called the strong generalized solution of equation (\ref{Pok-hyperbolic-eq:1}).
In such a way we obtain an operator $\overline{Q} = \overline{\mathcal{Q}}\mathcal{P}_{\overline{X}} : \overline{\mathcal{H}} \rightarrow \mathcal{H}_{1}$,
which is normally resolvable ($R\left(\overline{Q}\right) = \overline{R\left(\overline{Q}\right)}$) and has Moore--Penrose pseudoinverse $\overline{Q}^{+}$, $\overline{\mathcal{H}} = N(Q) \oplus \overline{X}$.
Then the set of strong generalized solutions of the equation (\ref{Pok-hyperbolic-eq:5}) has the form:
\begin{equation*}
c = \overline{Q}^{+}g_{1} + \mathcal{P}_{N\left(\overline{Q}\right)}\overline{c}, \quad \forall \overline{c} \in \overline{\mathcal{H}}.
\end{equation*}

Finally, consider the case when $g_{1} \notin \overline{R(Q)}$. This condition is equivalent to the following $\mathcal{P}_{N\left(\overline{Q}^{*}\right)}g_{1} \neq 0$.
In this case, there are elements $c\in\overline{\mathcal{H}}$
\begin{equation*}
c = \overline{Q}^{+}g_{1} + \mathcal{P}_{N\left(\overline{Q}\right)}\overline{c}, \quad \overline{c} \in
\overline{\mathcal{H}},
\end{equation*}
that minimize the norm  $||\overline{Q}c - g_{1} ||_{\overline{\mathcal{H}}}$. These elements are called strong pseudosolutions by analogy with \cite{Boi-Sam-2016}.

For the boundary-value problem (\ref{Pok-hyperbolic-eq:3}) the following assertion is true \cite{Boi-Pok-2019}.
\begin{thm}\label{Pok-hyperbolic-th-1}
1. a) The boundary-value problem (\ref{Pok-hyperbolic-eq:3}) possesses a strong generalized solutions if and only if the following condition is satisfied:
\begin{equation} \label{Pok-hyperbolic-eq:6}
\mathcal{P}_{\overline{\mathcal{H}}_{\overline{Q}}}g_{1} = 0;
\end{equation}
if
\begin{equation*}
g_{1} \in R(Q),
\end{equation*}
then strong generalized solutions are classical;

b) if the condition (\ref{Pok-hyperbolic-eq:6}) is satisfied, then the set of strong generalized solutions of boundary-value problem (\ref{Pok-hyperbolic-eq:3}) has the form:
\begin{equation} \label{Pok-hyperbolic-eq:7t}
x_{0}(t,c) = U(t)\mathcal{P}_{N\left(\overline{Q}\right)}c + \overline{(G[g, \alpha])(t)} \quad \forall c \in \overline{\mathcal{H}},
\end{equation}
where
\begin{equation*}
\overline{(G[g, \alpha])(t)} = U(t)\overline{Q}^{+}g_{1} + \int\limits_{0}^{t}U(t)U^{-1}(\tau)g(\tau)d\tau
\end{equation*}
is a generalized Green operator of the boundary-value problem (\ref{Pok-hyperbolic-eq:3}) which extended to space $\overline{\mathcal{H}}$;

2. a) The boundary-value problem (\ref{Pok-hyperbolic-eq:3}) possesses a strong generalized pseudosolutions if and only if the following condition is satisfied:
\begin{equation} \label{Pok-hyperbolic-eq:7}
\mathcal{P}_{\overline{\mathcal{H}}_{\overline{Q}}}g_{1} \neq 0;
\end{equation}

b) if the condition (\ref{Pok-hyperbolic-eq:7}) is satisfied, then the set of strong generalized pseudosolutions of boundary-value problem (\ref{Pok-hyperbolic-eq:3}) has the form:
\begin{equation*}
x_{0}(t,c) = U(t)\mathcal{P}_{N\left(\overline{Q}\right)}c + \overline{(G[g, \alpha])(t)} \quad \forall c \in \mathcal{\overline{H}},
\end{equation*}
where the elements of $c\in\overline{\mathcal{H}}$ minimize the norm $||\overline{Q}c - g_{1} ||_{\overline{\mathcal{H}}}$, $\mathcal{P}_{\overline{\mathcal{H}}_{\overline{Q}}}$ is orthprojector onto subspace $\overline{\mathcal{H}}_{\overline{Q}},$ $ \overline{\mathcal{H}} = \overline{\mathcal{H}}_{\overline{Q}} \oplus R(\overline{Q})$.
\end{thm}

\section{Conditions for the existence of the solution of the boundary-value problem}We consider the questions of finding necessary and sufficient conditions of solvability the nonlinear boundary-value problem (\ref{Pok-hyperbolic-eq:123}), (\ref{Pok-hyperbolic-eq:124}). To do this, we rewrite it as follows:
\begin{align}
&x^{\prime}(t, \varepsilon) = B(t)x(t, \varepsilon) + \varepsilon H(t, x(t, \varepsilon)) + g(t), \label{Pok-hyperbolic-eqc:1} \\
&lx(\cdot, \varepsilon) = \alpha, \label{Pok-hyperbolic-eqc:2}
\end{align}
where
\begin{equation*}
x(t, \varepsilon) = {\rm col} \begin{pmatrix}
	x_1(t, \varepsilon), & x_2(t, \varepsilon)
\end{pmatrix}, \quad x_1(t, \varepsilon) = y(t, \varepsilon), \quad x_2(t, \varepsilon) = y^{\prime}(t, \varepsilon),
\end{equation*}
and the nonlinearity $H(t, x(t, \varepsilon))$ has the form
\begin{equation*}
H(t, x(t, \varepsilon)) = \begin{pmatrix}
	0 \\
Z(t, x_1(t, \varepsilon), x_2(t, \varepsilon))
\end{pmatrix}.
\end{equation*}
We seek the solution $x(t, \varepsilon)$ of boundary-value problem (\ref{Pok-hyperbolic-eqc:1}), (\ref{Pok-hyperbolic-eqc:2}), which for $\varepsilon = 0$ turns into the solution $x_{0}(t, c)$ of the generating boundary-value problem (\ref{Pok-hyperbolic-eq:4}), (\ref{Pok-hyperbolic-eq:04}).
According to the Theorem \ref{Pok-hyperbolic-th-1}, the boundary-value problem (\ref{Pok-hyperbolic-eqc:1}), (\ref{Pok-hyperbolic-eqc:2}) possesses a strong generalized solutions if and only if the following condition is satisfied:
\begin{equation*}
\mathcal{P}_{\overline{\mathcal{H}}_{\overline{Q}}} \left\{g_{1} - \varepsilon l\int\limits_{0}^{\cdot}U(\cdot)U^{-1}(\tau)H(\tau, x(\tau, \varepsilon))d\tau \right\} = 0.
\end{equation*}
Using the condition (\ref{Pok-hyperbolic-eq:6}) and the fact that the nonlinearity $H$ is continuous in the neighborhood of the generating solution, that is $H(\tau, x(\tau, \varepsilon)) \rightarrow H(\tau, x_{0}(\tau, c))$, when $\varepsilon \rightarrow 0$, we obtain the following condition for the solvability of the nonlinear boundary-value problem (\ref{Pok-hyperbolic-eqc:1}), (\ref{Pok-hyperbolic-eqc:2}).

\begin{thm}\label{Pok-hyperbolic-th-2} If the boundary-value problem (\ref{Pok-hyperbolic-eqc:1}), (\ref{Pok-hyperbolic-eqc:2}) has a strong generalized solutions $x(t, \varepsilon)$, which for  $\varepsilon = 0$ turn into one of the solutions $x_{0}\left(t, c^{0}\right)$ generating boundary-value problem (\ref{Pok-hyperbolic-eq:4}), (\ref{Pok-hyperbolic-eq:04}), then the vector $c = c^{0}$ must be a solution to the operator equation for the generating elements
\begin{equation} \label{Pok-hyperbolic-eqvat:1}
F(c) = \mathcal{P}_{\overline{\mathcal{H}}_{\overline{Q}}}l\int\limits_{0}^{\cdot}U(\cdot)U^{-1}(\tau)H(\tau, x_{0}(\tau, c))d\tau = 0.
\end{equation}
\end{thm}

To obtain a sufficient condition for the existence of a solution, we make the following change of variables of the boundary-value problem (\ref{Pok-hyperbolic-eqc:1}), (\ref{Pok-hyperbolic-eqc:2}):
\begin{equation}\label{Pok-hyperbolic-eqvat:01}
x(t, \varepsilon) = y(t, \varepsilon) + x_{0}\left(t, c^{0}\right),
\end{equation}
where element $c^{0}$ is a solution to the operator equation for the generating elements (\ref{Pok-hyperbolic-eqvat:1}). By the change of variables (\ref{Pok-hyperbolic-eqvat:01}) in (\ref{Pok-hyperbolic-eqc:1}), we obtain the following boundary-value problem
\begin{align}
&y^{\prime}(t, \varepsilon) = B(t)y(t, \varepsilon) + \varepsilon H\left(t, y(t, \varepsilon) + x_{0}\left(t, c^{0}\right)\right), \label{Pok-hyperbolic-eqcu:1} \\
&ly(\cdot, \varepsilon) = 0. \label{Pok-hyperbolic-eqcu:2}
\end{align}

Decomposing the nonlinearity $H$ in the neighborhood of the generating solution
\begin{equation*}
H\left(t, y(t, \varepsilon) + x_{0}\left(t, c^{0}\right)\right) = H\left(t, x_{0}\left(t, c^{0}\right)\right) + H^{\prime}_{x}\left(t, x_{0}\left(t, c^{0}\right)\right)y(t, \varepsilon) + R(t, y(t, \varepsilon)),
\end{equation*}
where
\begin{equation*}
R(t, 0) = R^{\prime}_x(t, 0) = 0,
\end{equation*}
we arrive at the following boundary-value problem
\begin{align}
y^{\prime}(t, \varepsilon) = {}& B(t)y(t, \varepsilon) \notag\\
 {}&+ \varepsilon \left(H\left(t, x_{0}\left(t, c^{0}\right)\right) +  H^{\prime}_{x}\left(t, x_{0}\left(t, c^{0}\right)\right)y(t, \varepsilon) + R(t, y(t, \varepsilon))\right),  \label{Pok-hyperbolic-eqcut:1} \\
ly(\cdot, \varepsilon) = {}& 0. \label{Pok-hyperbolic-eqcut:2}
\end{align}

Necessary and sufficient condition for the existence of strong generalized solutions of the boundary-value problem (\ref{Pok-hyperbolic-eqcut:1}), (\ref{Pok-hyperbolic-eqcut:2}) is as follows:
\begin{align*}
\mathcal{P}_{\overline{\mathcal{H}}_{\overline{Q}}} {}& l\int\limits_{0}^{\cdot}U(\cdot)U^{-1}(\tau)\left( H\left(\tau, x_{0}\left(\tau, c^{0}\right)\right) \right.\\
 {}&\left. + H_{x}^{\prime}\left(\tau, x_{0}\left(\tau, c^{0}\right)\right)y(\tau, \varepsilon) + R(\tau, y(\tau, \varepsilon)))\right)d\tau = 0.
\end{align*}
According to (\ref{Pok-hyperbolic-eqvat:1}), this condition is equivalent to the following condition:
\begin{equation}  \label{Pok-hyperbolic-eqcut:4}
\mathcal{P}_{\overline{\mathcal{H}}_{\overline{Q}}}l\int\limits_{0}^{\cdot}U(\cdot)U^{-1}(\tau)\left( H_{x}^{\prime}\left(\tau, x_{0}\left(\tau, c^{0}\right)\right)y(\tau, \varepsilon) + R(\tau, y(\tau, \varepsilon)))\right)d\tau = 0.
\end{equation}
If the condition (\ref{Pok-hyperbolic-eqcut:4}) is satisfied, then the set of strong generalized solutions of boundary-value problem (\ref{Pok-hyperbolic-eqcut:1}), (\ref{Pok-hyperbolic-eqcut:2}) has the form:
\begin{equation} \label{Pok-hyperbolic-eqcut:5}
y(t, \varepsilon) = U(t)\mathcal{P}_{N\left(\overline{Q}\right)}c + \varepsilon \overline{G[H(\cdot, y + x_{0}),0]}(t).
\end{equation}
Substituting the solution (\ref{Pok-hyperbolic-eqcut:5}) in the condition of solvability (\ref{Pok-hyperbolic-eqcut:4}), we obtain the following
operator equation for $c$:
\begin{equation} \label{Pok-hyperbolic-eqcut:6}
B_{0}c = b,
\end{equation}
where an operator $B_{0}$ and an element $b$ have the form
\begin{equation*}
B_{0} = \mathcal{P}_{\overline{\mathcal{H}}_{\overline{Q}}}l\int\limits_{0}^{\cdot}U(\cdot)U^{-1}(\tau)H_{x}^{\prime}\left(\tau, x_{0}\left(\tau, c^{0}\right)\right)U(\tau)\mathcal{P}_{N\left(\overline{Q}\right)}d\tau,
\end{equation*}
\begin{align*}
b = {}& -\mathcal{P}_{\overline{\mathcal{H}}_{\overline{Q}}}l\int\limits_{0}^{\cdot} U(\cdot)U^{-1}(\tau)\left(R(\tau, y(\tau, \varepsilon)) \right.\\
 {}&\left. + \varepsilon H_{x}^{\prime}\left(\tau, x_{0}\left(\tau, c^{0}\right)\right) \overline{G[H(\cdot, y + x_{0}),0]}(\tau)\right)d\tau.
\end{align*}

Necessary and sufficient condition for the existence of solutions of the equation (\ref{Pok-hyperbolic-eqcut:6}) is as follows:
\begin{equation} \label{Pok-hyperbolic-eqcut:7ns}
\mathcal{P}_{\mathcal{H}_{\overline{Q}}}b = 0.
\end{equation}
If
\begin{equation} \label{Pok-hyperbolic-eqcut:7}
\mathcal{P}_{\overline{\mathcal{H}}_{\overline{B}_{0}}}\mathcal{P}_{\overline{\mathcal{H}}_{\overline{Q}}} = 0,
\end{equation}
then equality (\ref{Pok-hyperbolic-eqcut:7ns}) is always satisfied and operator equation
(\ref{Pok-hyperbolic-eqcut:6}) possesses a strong generalized solution
\begin{equation} \label{Pok-hyperbolic-eqcut:07}
c = \overline{B}_{0}^{+}b.
\end{equation}
Here $\mathcal{P}_{\overline{\mathcal{H}}_{\overline{B}_{0}}}$ is the orthoprojector onto the subspace $\overline{\mathcal{H}}_{\overline{B}_{0}}$ ($\overline{\mathcal{H}}_{\overline{B}_{0}} = \widetilde{\overline{\mathcal{H}}}_{\overline{Q}} \ominus R(\overline{B}_{0}), \widetilde{\overline{\mathcal{H}}}_{\overline{Q}}$ is a completion of $\overline{\mathcal{H}}_{\overline{Q}}$ by the corresponding norm). 
\begin{rmk}\label{Pok-hyperbolic-rem-1} If the condition (\ref{Pok-hyperbolic-eqcut:7ns}) is not satisfied, then the representation (\ref{Pok-hyperbolic-eqcut:07}) determines a strong generalized pseudosolution of the equation (\ref{Pok-hyperbolic-eqcut:6}), which minimize the norm $||\overline{B}_{0}c - b||_{\overline{\mathcal{H}}_{\overline{B}_{0}}}$.
\end{rmk}

Thus, we can rewrite the boundary-value problem (\ref{Pok-hyperbolic-eqcut:1}), (\ref{Pok-hyperbolic-eqcut:2}) as the following operator system
\begin{align*}
&y(t, \varepsilon) = U(t)\mathcal{P}_{N\left(\overline{Q}\right)}c + \overline{y}(t, \varepsilon), \\
&c = - \overline{B}_{0}^{+}\mathcal{P}_{\overline{\mathcal{H}}_{\overline{Q}}}l\int\limits_{0}^{\cdot} U(\cdot)U^{-1}(\tau) \left(H_{x}^{\prime}\left(\tau, x_{0}\left(\tau, c^{0}\right)\right) \overline{y}(\tau, \varepsilon) + R(\tau, y(\tau, \varepsilon))\right)d\tau, \\
&\overline{y}(t, \varepsilon) = \varepsilon \overline{G[H(\cdot, y + x_{0}),0]}(t),
\end{align*}
which can be solved using a convergent iterative process explained in detail in \cite{Boi-Sam-2016}. Hence, the following Theorem is true.

\begin{thm}\label{Pok-hyperbolic-th-3} Let the generating problem for (\ref{Pok-hyperbolic-eq:123}), (\ref{Pok-hyperbolic-eq:124}) problem (\ref{Pok-hyperbolic-eq:4}), (\ref{Pok-hyperbolic-eq:04}), subject to the condition (\ref{Pok-hyperbolic-eq:6}), have set of solutions $x_{0}(t,c)$ (\ref{Pok-hyperbolic-eq:7t}). Then for every real value of the vector $c=c^0$, which satisfies the equation for the generating elements (\ref{Pok-hyperbolic-eqvat:1}) and when conditions (\ref{Pok-hyperbolic-eqcut:7}) is satisfied, the boundary-value problem (\ref{Pok-hyperbolic-eq:123}), (\ref{Pok-hyperbolic-eq:124}) has a strong generalized solution $x(t,\varepsilon)$, which for $\varepsilon = 0$ turns into one of the solutions $x_{0}\left(t, c^{0}\right)$ generating boundary-value problem (\ref{Pok-hyperbolic-eq:4}), (\ref{Pok-hyperbolic-eq:04}). This solution can be found using the following iterative process:
\begin{align*}
&y_{k + 1}(t, \varepsilon) = U(t)\mathcal{P}_{N\left(\overline{Q}\right)}c_{k} + \overline{y}_{k+1}(t, \varepsilon), \\
&c_{k} = - \overline{B}_{0}^{+}\mathcal{P}_{\overline{\mathcal{H}}_{\overline{Q}}}l\int\limits_{0}^{\cdot}U(\cdot)U^{-1}(\tau) \left(H_{x}^{\prime}\left(\tau, x_{0}\left(\tau, c^{0}\right)\right) \overline{y}_{k}(\tau, \varepsilon) + R(\tau, y_{k}(\tau, \varepsilon))\right)d\tau, \\
&\overline{y}_{k + 1}(t, \varepsilon) = \varepsilon \overline{G[H(\cdot, y_{k} + x_{0}),0]}(t),
\end{align*}
\begin{equation*}
x_{k}(t, \varepsilon) = y_{k}(t, \varepsilon) + x_{0}\left(t, c^{0}\right), \quad x(t, \varepsilon) = \lim_{k \rightarrow \infty} x_{k}(t, \varepsilon), \quad y_{0}(t,\varepsilon) = \overline{y}_{0}(t,\varepsilon) = 0,
\end{equation*}
\begin{align*}
R(t, y_{k}(t, \varepsilon)) = {}& H\left(t, y_{k}(t, \varepsilon) + x_{0}\left(t, c^{0}\right)\right) \\
 {}& - H\left(t, x_{0}\left(t, c^{0}\right)\right) - H^{\prime}_{x}\left(t, x_{0}\left(t, c^{0}\right)\right)y_{k}(t, \varepsilon).
\end{align*}
\end{thm}

\section{Applications}
In this section, we give two test problems for to illustrate the theoretical results presented above: linear hyperbolic equation with constant operator and Van der Pol equation.

\subsection{Linear hyperbolic equation with constant operator}

We consider the evolutionary differential equation in the separable Hilbert space $\mathcal{H}$:
\begin{equation} \label{Pok-hyperbolic-equr:1}
y^{\prime\prime}(t) + Ty(t) = f(t),
\end{equation}
with boundary condition
\begin{equation} \label{Pok-hyperbolic-equr:2}
l_{1}(y(\cdot)) = \alpha_1, \quad l_{2}(y^{\prime}(\cdot)) = \alpha_{2},
\end{equation}
where $y \in C^2(J,\mathcal{H})$, $J=[0,w]\subset\mathbb{R}$, $T$ is an unbounded constant operator with compact inverse
$T^{-1}$, $f \in C(J,\mathcal{H})$, $l_{1}:C^2[J,\mathcal{H}]\rightarrow\mathcal{H}$, $l_{2}:C^1[J,\mathcal{H}]\rightarrow\mathcal{H}$ are bounded linear functionals, $\alpha_{1} \in \mathcal{H}_{1}$, $\alpha_{2} \in \mathcal{H}_{1}$. Hence, in this case we have
\begin{equation*}
l = \begin{pmatrix} l_1 & 0 \\
0 & l_2 \end{pmatrix}, \quad
\alpha = \begin{pmatrix} \alpha_1 \\ \alpha_2 \end{pmatrix}
\end{equation*}

We consider the case of the separable Hilbert space. Therefore, there exists the orthonormal basis $\{e_{i}\}_{i=1}^{\infty} \subset \mathcal{H}$ such that
\begin{equation*}
	y(t) = \sum_{i = 1}^{\infty}c_{i}(t)e_{i}, \quad Ty(t) = \sum_{i =
1}^{\infty}\lambda_{i}c_{i}(t)e_{i}, \quad \lambda_{i} \rightarrow
\infty,
\end{equation*}
\begin{equation*}
	f(t) = \sum_{i = 1}^{\infty}f_{i}(t)e_{i}, \quad c_{i}(t) = (y(t), e_{i}), \quad f_{i}(t) = (f(t), e_{i}).
\end{equation*}
The first of the boundary conditions (\ref{Pok-hyperbolic-equr:2}) can be rewritten as follows:
\begin{equation*}
l_1(y(\cdot)) = l_1\left(\sum_{i = 1}^{\infty}c_{i}(\cdot)e_{i}\right)  = \sum_{i = 1}^{\infty}l_1\left( c_{i}(\cdot)e_{i}\right).
\end{equation*}
According to $l(c_{i}(\cdot)e_{i}) \subset \mathcal{H}$, $i=\overline{1,\infty}$, the following expansions are true:
\begin{equation*}
l_1(c_{i}(\cdot)e_{i}) = \sum\limits_{j = 1}^{\infty}a_{ij}e_{j}, \quad i=\overline{1,\infty},
\end{equation*}
where
\begin{equation*}
a_{ij} = (l_1(c_{i}(\cdot)e_{i}), e_{j}), \quad \sum\limits_{i=1}^{\infty}a_{ij}^{2} < \infty, \quad j=\overline{1,\infty}.
\end{equation*}
Decomposing the element $\alpha_1$ on the basis of the space $\mathcal{H}$, we get
\begin{equation*}
\alpha_1 = \sum\limits_{j = 1}^{\infty} \alpha_{1 j}e_{j}.
\end{equation*}
Thus, the first boundary condition (\ref{Pok-hyperbolic-equr:2}) is equivalent to the countable number of conditions of the following form:
\begin{equation*}
\sum_{i = 1}^{\infty}a_{ij} = \alpha_{1j}, \quad j = \overline{1, \infty}.
\end{equation*}
Similarly, the second boundary condition (\ref{Pok-hyperbolic-equr:2}) is equivalent to the countable number of the following conditions:
\begin{equation*}
\sum_{i = 1}^{\infty}b_{ij} = \alpha_{2j}, \quad b_{ij} = (l_2(c^{\prime}_{i}(\cdot)e_{i}), e_{j}), \quad j = \overline{1, \infty}.
\end{equation*}
By the change of variables
\begin{equation*}
c_{k}(t) = x_{k}(t), \quad c^{\prime}_{k}(t) = \sqrt{\lambda_k}y_{k}(t), \quad k = \overline{1, \infty},
\end{equation*}
boundary-value problem (\ref{Pok-hyperbolic-equr:1}), (\ref{Pok-hyperbolic-equr:2}) can be represented as the following countable system of differential equations
\begin{equation}\label{Pok-hyperbolic-equr:3}
\begin{split}
&x_{k}^{\prime}(t) = \sqrt{\lambda_k}y_{k}(t), \\
&y_{k}^{\prime}(t) = -\sqrt{\lambda_{k}}x_{k}(t) + \frac{f_{k}(t)}{\sqrt{\lambda_{k}}}, \quad k = \overline{1, \infty}
\end{split}
\end{equation}
with boundary conditions
\begin{equation} \label{Pok-hyperbolic-equr:4}
\sum_{i = 1}^{\infty}a_{ij} = \alpha_{1j}, \quad \sum_{i = 1}^{\infty}b_{ij} = \alpha_{2j}, \quad j = \overline{1, \infty}.
\end{equation}
We consider the resonant case when $\lambda_{k} = 4\pi^{2}k^{2}/w^{2}$ and for simplicity $w=2\pi$. Introducing the vector function
\begin{equation*}
x_{0}(t) = {\rm col} \begin{pmatrix} x_{1}(t), & y_{1}(t), & x_{k}(t), & y_{k}(t), & \cdots \end{pmatrix},
\end{equation*}
we rewrite the boundary-value problem (\ref{Pok-hyperbolic-equr:3}), (\ref{Pok-hyperbolic-equr:4}) in the form
\begin{equation} \label{Pok-hyperbolic-equr:5}
x_{0}^{\prime}(t) = Bx_{0}(t) + g(t),
\end{equation}
\begin{equation} \label{Pok-hyperbolic-equr:06}
\sum_{i = 1}^{\infty}a_{ij} = \alpha_{1j}, \quad \sum_{i = 1}^{\infty}b_{ij} = \alpha_{2j}, \quad j = \overline{1, \infty},
\end{equation}
where
\begin{equation}\label{Pok-hyperbolic-rrr:1}
B = {\rm diag} \begin{pmatrix}
B_{1} & B_{2} & \cdots & B_{k} & \cdots \end{pmatrix}, \quad
B_{k} = \begin{pmatrix}
0 & k \\
 -k & 0
\end{pmatrix},
\end{equation}
\begin{equation*}
g(t) = {\rm col} \begin{pmatrix} 0, & f_{1}(t), & 0, & f_{2}(t)/2, & 0, & f_{k}(t)/k, & \cdots \end{pmatrix}, \quad k = \overline{1, \infty}.
\end{equation*}
The evolutionary operator of system (\ref{Pok-hyperbolic-equr:5}) has the following form:
\begin{equation*}
U(t) = {\rm diag} \begin{pmatrix}
U_{1}(t) & U_{2}(t) & \cdots & U_{k}(t) & \cdots
\end{pmatrix}, \quad
U_{k}(t) = \begin{pmatrix}
\cos kt & \sin kt \\
-\sin kt & \cos kt
\end{pmatrix}.
\end{equation*}
The general solution of the linear inhomogeneous system (\ref{Pok-hyperbolic-equr:5}) has the following form:
\begin{align}
x_{0}(t, c) = {}& U(t)c+ \int\limits_{0}^{t}U(t)U^{-1}(\tau)g(\tau)d\tau \notag\\
= {}& {\rm col} \begin{pmatrix}
x_{1}(t, c), & y_{1}(t, c), & \cdots, & x_{k}(t, c), & y_{k}(t, c), & \cdots
\end{pmatrix}, \label{Pok-hyperbolic-eqmat:1}
\end{align}
where
\begin{equation*}
\begin{split}
&x_{k}(t, c)= c_{2k - 1}\cos k t + c_{2k}\sin k t + \frac{1}{k}\int\limits_{0}^{t}\sin k(t - \tau)f_{k}(\tau)d\tau, \\
&y_{k}(t, c)= -c_{2k - 1}\sin k t + c_{2k}\cos k t + \frac{1}{k}\int\limits_{0}^{t}\cos k(t - \tau)f_{k}(\tau)d\tau, \quad k = \overline{1, \infty}.
\end{split}
\end{equation*}
Substituting (\ref{Pok-hyperbolic-eqmat:1}) in boundary conditions (\ref{Pok-hyperbolic-equr:06}), we obtain that the elements $a_{ij}$, $b_{ij}$, $i,j = \overline{1, \infty}$ have the form:
\begin{align*}
a_{ij} = {}& c_{2i - 1}l_1(\cos i(\cdot)e_i, e_j) + c_{2i}l_1(\sin i(\cdot)e_i, e_j) \\
 {}&  + \frac{1}{i} l_1\left(\int\limits_{0}^{\cdot}\sin i(\cdot - \tau)f_{i}(\tau))d\tau e_{i}, e_{j}\right), \\
b_{ij} = {}& -i c_{2i - 1}l_2(\sin i(\cdot)e_{i}, e_{j}) + i c_{2i}l_2(\cos i(\cdot)e_i, e_j) \\
 {}& + l_2\left(\int\limits_{0}^{\cdot}\cos i(\cdot - \tau)f_{i}(\tau)d\tau e_{i}, e_{j}\right).
\end{align*}
We are noting
\begin{equation*}
c_{1}^{i} = c_{2i - 1}, \quad c_{2}^{i} = c_{2i}, \quad
d_{(2i)j} = l_{1}(\sin i(\cdot)e_i, e_j), \quad d_{(2i - 1)j} = l_1(\cos i(\cdot)e_{i}, e_{j}),
\end{equation*}
\begin{equation*}
q_{(2i)j} = i l_2(\cos i(\cdot)e_i, e_j), \quad  g^{1}_{1j} = \alpha_{1j} - \sum\limits_{i = 1}^{\infty}\frac{1}{i}l_{1}\left(\int\limits_{0}^{\cdot}\sin i(\cdot - \tau)f_{i}(\tau)d\tau e_{i}, e_{j}\right),
\end{equation*}
\begin{equation*}
 q_{(2i - 1) j} = -i l_2(\sin i(\cdot)e_i, e_j), \quad g^{1}_{2j} = \alpha_{2j} - \sum\limits_{i = 1}^{\infty}l_2\left(\int\limits_{0}^{\cdot}\cos i(\cdot - \tau)f_{i}(\tau)d\tau e_i, e_{j}\right).
\end{equation*}

Thus, finding the constants $c_{1}^{i}$, $c_{2}^{i}$, which are determined by boundary conditions (\ref{Pok-hyperbolic-equr:06}) is equivalent to solving the such countable system of linear of algebraic equations in the space of infinite sequences $\ell_{2}$:
\begin{equation} \label{Pok-hyperbolic-equr:6}
Qc = g_{1},
\end{equation}
where
\begin{equation*}
Q =  \begin{pmatrix}
Q_{11} & Q_{12} & \cdots & Q_{1k} & \cdots  \\
Q_{21} & Q_{22} & \cdots & Q_{2k} & \cdots  \\
\vdots & \vdots  & \ddots & \vdots & \vdots  \\
Q_{k1} & Q_{k2} & \cdots & Q_{kk} & \cdots  \\
\cdots  & \cdots  & \cdots & \cdots  & \cdots  \
\end{pmatrix}, \quad
Q_{ij} = \begin{pmatrix}
d_{(2i - 1)j} & d_{(2i)j} \\
q_{(2i - 1)j} & q_{(2i)j} \
\end{pmatrix},
\end{equation*}
\begin{equation*}
\begin{split}
c = {}& {\rm col} \begin{pmatrix} c_{1}^{1}, & c_{2}^{1}, & c_{1}^{2}, & c_{2}^{2}, & \cdots, & c_{1}^{k}, & c_{2}^{k}, & \cdots \end{pmatrix}, \\
g_{1} = {}& {\rm col} \begin{pmatrix} g^{1}_{11}, & g^{1}_{21}, & g^{1}_{12}, & g^{1}_{22}, & \cdots, & g^{1}_{1k}, & g^{1}_{2k}, & \cdots \end{pmatrix}.
\end{split}
\end{equation*}
If the set of values of the operator $Q$, which is determined by the matrix in the equation (\ref{Pok-hyperbolic-equr:6}), is closed $R(Q) = \overline{R(Q)}$, then the necessary and sufficient  solvability condition of the equation (\ref{Pok-hyperbolic-equr:6}) is as follows:
\begin{equation} \label{Pok-hyperbolic-equr:7}
\mathcal{P}_{\mathcal{H}_{Q}}g_{1} = 0,
\end{equation}
where $\mathcal{P}_{\mathcal{H}_{Q}} = I - QQ^{+}$ is the orthoprojector onto the kernel of the operator $Q^{*}$, $Q^{+}$ is the pseudoinverse Moore--Penrose matrix. Then the set of solutions of the equation (\ref{Pok-hyperbolic-equr:6}) has the form:
\begin{equation} \label{Pok-hyperbolic-equr:8}
c =Q^{+}g_{1}+\mathcal{P}_{N(Q)}h \quad \forall h \in \ell_{2},
\end{equation}
where $\mathcal{P}_{N(Q)} = I - Q^{+}Q$ is the orthoprojector onto the kernel of the operator $Q$.
Substituting the vector (\ref{Pok-hyperbolic-equr:8}) in the representation (\ref{Pok-hyperbolic-eqmat:1}), we obtain a general classical generalized solution of the original boundary-value problem (\ref{Pok-hyperbolic-equr:1}), (\ref{Pok-hyperbolic-equr:2}):
\begin{equation} \label{Pok-hyperbolic-equr:9}
x_{0}(t, h) = U(t)\mathcal{P}_{N(Q)}h + (G[f, \alpha_1, \alpha_2])(t),
\end{equation}
where
\begin{equation} \label{Pok-hyperbolic-equr:10}
(G[f, \alpha_{1}, \alpha_2])(t) = U(t)Q^{+}g_1 + \int\limits_{0}^{t}
\begin{pmatrix}\sin(t - \tau)f_1(\tau) \\
\cos(t-\tau)f_1(\tau) \\
\sin 2(t - \tau)f_2(\tau)/2 \\
\cos 2(t - \tau)f_2(\tau)/2 \\
\cdots \\
\sin k(t - \tau)f_{k}(\tau)/k\\
\cos k(t - \tau)f_{k}(\tau)/k\\
\cdots
\end{pmatrix} d\tau
\end{equation}
is the generalized Green operator of the corresponding boundary-value problem.

If the set of values of the operator $Q$ is not closed, then there is the strong generalized operator  $\overline{Q}^{+}$ and the necessary and sufficient condition of the strong generalized solvability of the equation (\ref{Pok-hyperbolic-equr:6}) is as follows:
\begin{equation} \label{Pok-hyperbolic-equr:11}
\mathcal{P}_{N\left(\overline{Q}^{*}\right)}g_{1} = 0.
\end{equation}

If $g_{1} \in R(Q)$, then we obtain the classical generalized solution as in the previous case. If the condition of solvability (\ref{Pok-hyperbolic-equr:11}) is satisfied, then the set of strong generalized solutions of the equation (\ref{Pok-hyperbolic-equr:6}) has the form:
\begin{equation} \label{Pok-hyperbolic-equr:12}
c = \overline{Q}^{+}g_{1} + \mathcal{P}_{N\left(\overline{Q}\right)}h \quad, ~~ \in \ell_2.
\end{equation}

Substituting the vector (\ref{Pok-hyperbolic-equr:12}) in (\ref{Pok-hyperbolic-eqmat:1}), we obtain a general strong generalized solution of the initial boundary-value problem (\ref{Pok-hyperbolic-equr:1}), (\ref{Pok-hyperbolic-equr:2}):
\begin{equation} \label{Pok-hyperbolic-equr:09}
x_{0}(t, h) = U(t)\mathcal{P}_{N\left(\overline{Q}\right)}h + \overline{(G[f,\alpha_1, \alpha_2])}(t),
\end{equation}
where
\begin{equation} \label{Pok-hyperbolic-equr:010}
\overline{(G[f, \alpha_1, \alpha_2])}(t) = U(t)\overline{Q}^{+}g_{1} + \int\limits_{0}^{t} \begin{pmatrix}\sin(t - \tau)f_1(\tau) \\
\cos(t-\tau)f_1(\tau) \\
\sin 2(t - \tau)f_2(\tau)/2 \\
\cos 2(t - \tau)f_2(\tau)/2 \\
\cdots \\
\sin k(t - \tau)f_{k}(\tau)/k \\
\cos k(t - \tau)f_{k}(\tau)/k \\
\cdots\
\end{pmatrix} d\tau
\end{equation}
is the generalized Green operator of the corresponding boundary-value problem.

If condition
\begin{equation} \label{Pok-hyperbolic-equr:011}
\mathcal{P}_{\overline{\mathcal{H}}_{\overline{Q}}}g_{1} \neq 0,
\end{equation}
is satisfied, then the expression (\ref{Pok-hyperbolic-equr:09}) describes strong generalized pseudo-solutions of the boundary-value problem (\ref{Pok-hyperbolic-equr:1}), (\ref{Pok-hyperbolic-equr:2}).

Thus, we have proved the following theorem.
\begin{thm}\label{Pok-hyperbolic-th-4}
1 a. The boundary-value problem (\ref{Pok-hyperbolic-equr:1}), (\ref{Pok-hyperbolic-equr:2}) possesses a strong generalized solutions if and only if the condition (\ref{Pok-hyperbolic-equr:11}) is satisfied; if $g_{1} \in R(Q)$, then the obtained solutions are classical generalized solutions.

1 b. If the condition (\ref{Pok-hyperbolic-equr:11}) is satisfied, then the strong generalized solutions have the form (\ref{Pok-hyperbolic-equr:09}).

2 a. The boundary-value problem (\ref{Pok-hyperbolic-equr:1}), (\ref{Pok-hyperbolic-equr:2}) possesses a strong generalized pseudosolutions if and only if the condition (\ref{Pok-hyperbolic-equr:011}) is satisfied.

2 b. If the condition (\ref{Pok-hyperbolic-equr:011}) is satisfied, then the strong generalized pseudosolutions have the form (\ref{Pok-hyperbolic-equr:09}).
\end{thm}
Thus, the boundary-value problem (\ref{Pok-hyperbolic-equr:1}), (\ref{Pok-hyperbolic-equr:2}) is always solvable in one of the above senses.

\subsection{Van der Pol equation in Hilbert space}

We consider the evolutionary differential equation in the separable Hilbert space $\mathcal{H}$:
 \begin{equation} \label{Pok-hyperbolic-eqde:11}
 y^{\prime\prime}(t,\varepsilon) + T y(t, \varepsilon) = \varepsilon \left(1 -
||y(t, \varepsilon)||^{2}\right)y^{\prime}(t, \varepsilon),
\end{equation}
\begin{equation} \label{Pok-hyperbolic-eqde:12}
y(0, \varepsilon) = y(w, \varepsilon), \quad y^{\prime}(0, \varepsilon) = y^{\prime}(w, \varepsilon),
\end{equation}
where $y \in C^2(J,\mathcal{H})$, $J=[0,w]\subset\mathbb{R}$, $T$ is an unbounded constant operator with compact
inverse $T^{-1}$, $\varepsilon<<1$ is a small parameter.

In our case,
\begin{equation*}
f(t)=0, \quad Z(t,y(t,\varepsilon), y^{\prime}(t,\varepsilon))=\left(1 -
||y(t, \varepsilon)||^{2}\right)y'(t, \varepsilon),
\end{equation*}
\begin{equation*}
l(y(\cdot, \varepsilon), y^{\prime}(\cdot, \varepsilon)) =  \begin{pmatrix} y(0, \varepsilon) - y(w, \varepsilon) \\ y^{\prime}(0, \varepsilon) - y^{\prime}(w, \varepsilon) \end{pmatrix} = \begin{pmatrix} 0 \\ 0 \end{pmatrix} = \alpha.
\end{equation*}

We consider the case of the separable Hilbert space. Then there is an orthonormal basis $\{e_{i}\}_{i=1}^{\infty} \subset \mathcal{H}$ such that
\begin{equation*}
	y(t) = \sum_{i = 1}^{\infty}c_{i}(t)e_{i}, \quad Ty(t) = \sum_{i =
1}^{\infty}\lambda_{i}c_{i}(t)e_{i}, \quad \lambda_{i} \rightarrow
\infty.
\end{equation*}
The operator system (\ref{Pok-hyperbolic-eqcut:1}), (\ref{Pok-hyperbolic-eqcut:2}) for the boundary-value problem (\ref{Pok-hyperbolic-eqde:11}), (\ref{Pok-hyperbolic-eqde:12}) in this case, by the change of variables
\begin{equation*}
c_{k}(t) = x_{k}(t), \quad c^{\prime}_{k}(t) = \sqrt{\lambda_k}y_{k}(t), \quad k = \overline{1, \infty},
\end{equation*}
equivalent to a countable system of ordinary differential
equations of the form:
\begin{equation}\label{Pok-hyperbolic-eqde:14}
\begin{split}
&x^{\prime}_{k}(t) = \sqrt{\lambda_{k}}y_{k}(t), \\
&y^{\prime}_{k}(t) = -\sqrt{\lambda_{k}}x_{k}(t) +
 \frac{\varepsilon}{\sqrt{\lambda_{k}}}\left(1 - \sum_{j =
1}^{\infty}x_{j}^{2}(t)\right)y_{k}(t), \quad k = \overline{1, \infty}
\end{split}
\end{equation}
with boundary conditions
\begin{equation}\label{Pok-hyperbolic-eqde:15}
x_{k}(0) = x_{k}(w), \quad y_{k}(0)=y_{k}(w).
\end{equation}
We will find the solutions of this boundary-value problem in the space $C^{1}([0; w])$, for $\varepsilon = 0$ turn into one of the solutions generating boundary-value problem
\begin{equation}\label{Pok-hyperbolic-eqde1:11}
\begin{split}
&x^{\prime}_{k}(t) = \sqrt{\lambda_{k}}y_{k}(t), \\
&y^{\prime}_{k}(t) = -\sqrt{\lambda_{k}}x_{k}(t), \quad k = \overline{1, \infty},
\end{split}
\end{equation}
\begin{equation} \label{Pok-hyperbolic-eqde1:12}
x_{k}(0) = x_{k}(w), \quad y_{k}(0)=y_{k}(w).
\end{equation}

Consider a critical case $\lambda_{k} = 4\pi^{2}k^{2}/w^{2}$, $k \in \mathbb{N}$. Let, for simplicity, $w=2\pi$. In this case, the set of periodic solutions of the generating boundary-value problem (\ref{Pok-hyperbolic-eqde:14}), (\ref{Pok-hyperbolic-eqde:15}) has the form:
\begin{equation*}
x_{k}^{0}(t) = c_{1}^{k}\cos kt + c_{2}^{k}\sin kt, \quad y_{k}^{0}(t) = -c_{1}^{k}\sin kt + c_{2}^{k}\cos kt,
\end{equation*}
for all pairs of constants $c_{1}^{k}$, $c_{2}^{k} \in \mathbb{R}$, $k \in \mathbb{N}$. In this case $Q = U(0) - U(2\pi) = 0$, where the evolutionary operator has the form (\ref{Pok-hyperbolic-rrr:1}). Therefore, $\mathcal{P}_{\mathcal{H}_{Q}} = I$ and the equation for generating amplitudes can be written as an infinite vector $F(c)$:
\begin{equation} \label{porodg:1}
F(c): = {\rm col} \begin{pmatrix} F_{1}^{1}(c), & F_{2}^{1}(c), & F_{1}^{2}(c), & F_{2}^{2}(c), & \cdots \end{pmatrix} = 0,
\end{equation}
consisting of pairs $F_{1}^{k}(c)$, $F_{2}^{k}(c)$, $k \in \mathbb{N}$
\begin{align*}
F_{1}^{k}(c) = {}& -\frac{1}{k} \int\limits_{0}^{2\pi} \sin k\tau \left(1 - \sum_{j = 1}^{\infty}\left(c_{1}^{j}\cos j\tau + c_{2}^{j}\sin j\tau\right)^{2}\right) \\
{}& \times \left(-c_{1}^{k}\sin k\tau + c_{2}^{k}\cos k\tau\right)d\tau,
\end{align*}
\begin{align*}
F_{2}^{k}(c) = {}& \frac{1}{k} \int\limits_{0}^{2\pi} \cos k\tau \left(1 - \sum_{j = 1}^{\infty}\left(c_{1}^{j}\cos j\tau + c_{2}^{j}\sin j\tau\right)^{2}\right) \\
{}& \times \left(-c_{1}^{k} \sin k\tau+ c_{2}^{k}\cos k\tau\right)d\tau.
\end{align*}
After the transformations we obtain that the equation for the generating amplitudes will be equivalent to such a countable system of nonlinear algebraic equations
\begin{equation}\label{porodg:01}
F_{1}^{k}(c) = -\frac{\pi c_{1}^{k}}{4k}\left(2 \sum_{j = 1}^{\infty}\left(\left(c_{1}^{j}\right)^{2} + \left(c_{2}^{j}\right)^{2}\right) -
\left(c_{1}^{k}\right)^{2} - \left(c_{2}^{k}\right)^{2} - 4\right)
= 0,
\end{equation}
\begin{equation}\label{porodg:001}
F_{2}^{k}(c) =
-\frac{\pi c_{2}^{k}}{4k}\left(2 \sum_{j = 1}^{\infty}\left(\left(c_{1}^{j}\right)^{2} + \left(c_{2}^{j}\right)^{2}\right) -
\left(c_{1}^{k}\right)^{2} - \left(c_{2}^{k}\right)^{2} - 4\right) = 0,
\end{equation}
According to (\ref{porodg:01}), we have that $c_{1}^{k} = 0$ or
\begin{equation}\label{Pok-hyperbolic-equat:01}
2 \sum_{j = 1}^{\infty}\left(\left(c_{1}^{j}\right)^{2} + \left(c_{2}^{j}\right)^{2}\right) - \left(c_{1}^{k}\right)^{2} -
\left(c_{2}^{k}\right)^{2} - 4 = 0, \quad k \in \mathbb{N}.
\end{equation}
Similarly, according to (\ref{porodg:001}), we have that $c_{2}^{k} = 0$ or equality (\ref{Pok-hyperbolic-equat:01}) is true.
Take an arbitrary $m \neq k$. Then, we have that $c_{1}^{m} = 0$, or $c_{2}^{m} = 0$, or
\begin{equation} \label{Pok-hyperbolic-equat:2}
2 \sum_{j = 1}^{\infty}\left(\left(c_{1}^{j}\right)^{2} + \left(c_{2}^{j}\right)^{2}\right) - \left(c_{1}^{m}\right)^{2} -
\left(c_{2}^{m}\right)^{2} - 4 = 0, \quad m \in \mathbb{N}.
\end{equation}
Considering the difference (\ref{Pok-hyperbolic-equat:01}) and (\ref{Pok-hyperbolic-equat:2}), we get
\begin{equation} \label{Pok-hyperbolic-equat:3}
\left(c_{1}^{m}\right)^{2} + \left(c_{2}^{m}\right)^{2} = \left(c_{1}^{k}\right)^{2} + \left(c_{2}^{k}\right)^{2}.
\end{equation}
Let $\left(c_{1}^{m}\right)^{2} + \left(c_{2}^{m}\right)^{2} = a^{2}$. Necessary condition for the convergence of the series with the left part (\ref{Pok-hyperbolic-equat:01}) is the assumption that only a finite number of components $(c_{1}^{k}, c_{2}^{k})$, $k\in\mathbb{N}$ will be nonzero (otherwise you will get an infinite series consisting of the same numbers). That is, the infinite-dimensional generating boundary-value problem (\ref{Pok-hyperbolic-eqde1:11}), (\ref{Pok-hyperbolic-eqde1:12}) has a $N$-parameter family of solutions. Let $(c_{1}^{k_{i}}, c_{2}^{k_{i}}) \neq (0,0)$, $i = \overline{1,N}$ (or one of the components of this pair is nonzero). Then, using (\ref{Pok-hyperbolic-equat:2}), (\ref{Pok-hyperbolic-equat:3}) we get
\begin{equation*}
(2N -1)a^{2} = 4.
\end{equation*}
From here we finally get
\begin{equation}\label{Pok-hyperbolic-eq:1000}
\left(c_{1}^{k_{i}}\right)^{2} + \left(c_{2}^{k_{i}}\right)^{2} =
\left(\frac{2}{\sqrt{2N-1}}\right)^{2}, \quad  i = \overline{1,N}.
\end{equation}
Note that the problem of periodic solutions of the Van der Pol equation makes it possible to say that these constants $(c_{1}^{k_{i}},
c_{2}^{k_{i}})$, $i = \overline{1,N}$ are the amplitudes of the periodic solutions of the original Van der Pol equation. Thus, we obtained the following result.

\begin{thm}\label{Pok-hyperbolic-th-5} {\rm\textbf{(Necessary condition for the solvability of the Van der Pol equation)}}.
Let the boundary-value problem (\ref{Pok-hyperbolic-eqde:14}), (\ref{Pok-hyperbolic-eqde:15}) has a solution, which for $\varepsilon = 0$ turns into one of the solutions generating boundary-value problem (\ref{Pok-hyperbolic-eqde1:11}), (\ref{Pok-hyperbolic-eqde1:12}) with a set of pairs $(c_{1}^{k}, c_{2}^{k})$, $k \in \mathbb{N}$. Then among these pairs may be no more than a finite number of nonzero. Moreover, if $(c_{1}^{k_{i}}, c_{2}^{k_{i}}) \neq (0, 0)$, $i = \overline{1,N}$, then these constants are on the N-dimensional torus of a finite-dimensional subspace of constants (\ref{Pok-hyperbolic-eq:1000}).
\end{thm}

To obtain a sufficient condition for the existence of a solution, we will use the results on the existence of periodic solutions of the linear equation (\ref{Pok-hyperbolic-eqde:11}), i.e.
\begin{equation*}
 y^{\prime\prime}(t,\varepsilon) + T y(t, \varepsilon) = f(t),
\end{equation*}
where $f(t)\neq 0$. In this case (it follows from the first example), the generating boundary-value problem (\ref{Pok-hyperbolic-eqde1:11}), (\ref{Pok-hyperbolic-eqde1:12}) is solvable if and only if the conditions
\begin{equation*}
\int\limits_{0}^{2\pi} \sin k\tau f_{k}(\tau)d\tau = 0, \quad \int\limits_{0}^{2\pi} \cos k\tau f_{k}(\tau)d\tau = 0, \quad k=\overline{1,\infty}.
\end{equation*}
are satisfied.
Then, the set of solutions of the generating boundary-value problem has the form
\begin{equation*}
\begin{split}
&x_{k}^{0}\left(t, c_{1}^{k}, c_{2}^{k}\right) = c_{1}^{k}\cos kt + c_{2}^{k}\sin kt + \frac{1}{k}\int\limits_{0}^{t}\sin k(t - \tau)f_{k}(\tau)d\tau, \\
&y_{k}^{0}\left(t, c_{1}^{k}, c_{2}^{k}\right) = -c_{1}^{k}\sin kt + c_{2}^{k}\cos kt + \frac{1}{k}\int\limits_{0}^{t}\cos k(t - \tau)f_{k}(\tau)d\tau, \quad k=\overline{1,\infty}
\end{split}
\end{equation*}
and we get results from the work \cite{Boi-Pok-2019}.

Make the change of variables in the abstract Van der Pol equation (\ref{Pok-hyperbolic-eqde:11}), (\ref{Pok-hyperbolic-eqde:12})
\begin{equation*}
\begin{split}
&x_{k}(t,\varepsilon) = x_{k}^{0}\left(t, c^{k0}\right) + x_{k}^{1}(t, \varepsilon), \\
&y_{k}(t, \varepsilon) = y_{k}^{0}\left(t, c^{k0}\right) + y_{k}^{1}(t, \varepsilon),
\end{split}
\end{equation*}
where $x_{k}^{0}$, $y_{k}^{0}$ is generating solution, which is the solution of the boundary-value problem (\ref{Pok-hyperbolic-equr:3}) ($f(t) = 0$), and the vector of constants $c^{k0}$ satisfies the conditions of the Theorem \ref{Pok-hyperbolic-th-5}. In what follows, we assume that the first block of constants $(c_{1}^{i0}, c_{2}^{i0})$, $i = \overline{1, N}$ is nonzero. Then we come to the following a countable system of differential equations:
\begin{equation} \label{ura:1}
z^{\prime}(t, \varepsilon) = Bz(t,\varepsilon) + \varepsilon H\left(z(t, \varepsilon) + z_0\left(t, c^{0}\right)\right)
\end{equation}
or in the form
\begin{align}
z^{\prime}(t, \varepsilon) = {}& Bz(t, \varepsilon) + \varepsilon H\left(z_0\left(t, c^{0}\right)\right) \notag\\
 {}& + \varepsilon A_{1}(t)z(t, \varepsilon) + \varepsilon A_{2}(t)z(t, \varepsilon) + \varepsilon R(z(t, \varepsilon)), \label{eqva:1}
\end{align}
with boundary condition
\begin{equation} \label{eqva:2}
z(0, \varepsilon) = z(2\pi, \varepsilon).
\end{equation}
Here, the matrix $B$ has the form (\ref{Pok-hyperbolic-rrr:1}) and the vector functions $z(t, \varepsilon)$, $z_0\left(t, c^{0}\right)$, $H\left(z(t, \varepsilon) + z_0\left(t, c^{0}\right)\right)$ have the form:
\begin{equation*}
z(t, \varepsilon) = {\rm col} \begin{pmatrix} x_{1}^{1}(t, \varepsilon), & y_{1}^{1}(t, \varepsilon), & \cdots, & x_{k}^{1}(t, \varepsilon), & y_{k}^{1}(t, \varepsilon), & \cdots \end{pmatrix},
\end{equation*}
\begin{equation*}
z_{0}\left(t, c^{0}\right) = {\rm col} \begin{pmatrix} x_{0}\left(t, c^{0}\right), & y_{0}\left(t, c^{0}\right) \end{pmatrix},
\end{equation*}
\begin{equation*}
H(z(t, \varepsilon)) = {\rm col} \begin{pmatrix} 0, & H_{1}(z(t, \varepsilon)), & \cdots, & 0, & H_{k}(z(t, \varepsilon)), & \cdots \end{pmatrix},
\end{equation*}
\begin{equation*}
H_{k}(z(t, \varepsilon)) = \frac{1}{k}\left(1 - \sum\limits_{j = 1}^{\infty}\left(x_{j}^{1}(t, \varepsilon)\right)^2\right)y_{k}^{1}(t, \varepsilon),
\end{equation*}
where components of the vector-function $H\left(z(t, \varepsilon) + z_0\left(t, c^{0}\right)\right)$ have the following form:
\begin{equation*}
H\left(z_0\left(t, c^{0}\right)\right) = {\rm col} \begin{pmatrix} 0 & H_{1}\left(z_0\left(t, c^{0}\right)\right) & \cdots, & 0, & H_{k}\left(z_0\left(t, c^{0}\right)\right), & \cdots \end{pmatrix},
\end{equation*}
\begin{equation*}
H_{k}\left(z_0\left(t, c^{0}\right)\right) = \frac{1}{k}\left(1 - \sum\limits_{j = 1}^{N}\left(x_{j}^{0}\left(t, c^{0}\right)\right)^2\right)y_{k}^{0}\left(t, c^{0}\right) =
\end{equation*}
\begin{equation*}
 = \frac{1}{k}\left(1 - \sum\limits_{j = 1}^{N}\left(c_{1}^{j0}\cos jt + c_{2}^{j0}\sin jt\right)^{2}\right)\left(-c_{1}^{k0}\sin t + c_{2}^{k0}\cos t\right),
\end{equation*}
\begin{equation*}
 A_{1}(t) = \begin{pmatrix}
A^{1}_{11}(t) & A^{1}_{12}(t) & \cdots & A^{1}_{1N}(t) & O & \cdots  \\
A^{1}_{21}(t) & A^{1}_{22}(t) & \cdots & A^{1}_{2N}(t) & O & \cdots  \\
\vdots & \vdots  & \ddots & \vdots & \vdots  & \vdots  \\
A^{1}_{k1}(t) & A^{1}_{k2}(t) & \cdots & A^{1}_{NN}(t) & O & \cdots  \\
O & O & \cdots & O & O & \cdots  \\
\cdots  & \cdots  & \cdots & \cdots  & \cdots & \cdots
\end{pmatrix},
\end{equation*}
\begin{equation*}
 A_{2}(t) =  {\rm diag} \begin{pmatrix}
A^{2}_{1}(t) & A^{2}_{2}(t) & \cdots & A^{2}_{N}(t) & O & \cdots
\end{pmatrix},
\end{equation*}
\begin{equation*}
A^{1}_{ij}(t) =\frac{2}{i}x_{j}^{0}\left(t, c^{0}\right)y_{i}^{0}\left(t, c^{0}\right)
\begin{pmatrix}
0 & 0 \\
1 & 0
\end{pmatrix},
\end{equation*}
\begin{equation*}
A^{2}_{k}(t) =\frac{1}{k}\left(1 - \sum\limits_{j = 1}^{N}\left(x_{j}^{0}\left(t, c^{0}\right)\right)^2\right)
\begin{pmatrix}
0 & 0 \\
0 & 1
\end{pmatrix},
\end{equation*}
\begin{equation*}
R(z(t, \varepsilon)) = {\rm col}  \begin{pmatrix} 0, & R_{1}(z(t, \varepsilon)), & \cdots, & 0, & R_{k}(z(t, \varepsilon)), & \cdots  \end{pmatrix},
\end{equation*}
\begin{align*}
R_{k}(z(t, \varepsilon)) = {}& -\frac{1}{k}\sum\limits_{j = 1}^{\infty}\left(x_{j}^{1}(t, \varepsilon)\right)^2\left(y_{k}^{0}\left(t, c^{0}\right)+y_{k}^{1}(t, \varepsilon)\right) \\
 {}& -\frac{1}{k} \sum\limits_{j = 1}^{N}x_{j}^{0}\left(t, c^{0}\right)x_{j}^{1}(t, \varepsilon)y_{k}^{1}(t, \varepsilon).
\end{align*}

In this case, the condition for the solvability of the boundary-value problem (\ref{eqva:1}), (\ref{eqva:2}) has the form
\begin{equation} \label{eqva:3}
\int\limits_{0}^{2\pi}U^{-1}(\tau)\left( A_{1}(\tau)z(\tau, \varepsilon) + A_2(\tau)z(\tau, \varepsilon) + R(z(\tau, \varepsilon))\right)d\tau = 0.
\end{equation}
According to (\ref{eqva:3}), the solution of the boundary-value problem (\ref{eqva:1}), (\ref{eqva:2}) has the form
\begin{equation} \label{eqva:4}
z(t, \varepsilon) = U(t)c + \overline{z}(t, \varepsilon),
\end{equation}
where
\begin{equation*}
\overline{z}(t, \varepsilon) = \varepsilon G\left[H\left(z(\cdot, \varepsilon) + z_{0}\left(\cdot, c^{0}\right)\right)\right](t),
\end{equation*}
or in the form
\begin{equation*}
\overline{z}(t, \varepsilon) = \varepsilon G\left[H\left(z_{0}\left(\cdot, c^{0}\right)\right) + \left(A_{1}(\cdot) + A_2(\cdot)\right)z(\cdot, \varepsilon) + R(z(\cdot, \varepsilon))\right](t).
\end{equation*}
In the expanded form we have
\begin{equation}\label{ura:2}
\overline{z}(t, \varepsilon) = \varepsilon \int\limits_{0}^{t}\overline{H}(t, \tau, \varepsilon)d\tau
\end{equation}
and expression $\overline{H}(t, \tau, \varepsilon) = U(t - \tau)H\left(z(\tau, \varepsilon) + z_{0}\left(\tau, c^{0}\right)\right)$ is equal to
\begin{equation*}
\overline{H}(t, \tau, \varepsilon) = {\rm col} \begin{pmatrix} \overline{H}_{1}(t, \tau, \varepsilon),& \overline{H}_{2}(t, \tau, \varepsilon), & \cdots, & \overline{H}_{2k-1}(t, \tau, \varepsilon), & \overline{H}_{2k}(t, \tau, \varepsilon), & \cdots \end{pmatrix},
\end{equation*}
\begin{equation*}
\overline{H}_{2k-1}(t, \tau, \varepsilon) = \frac{1}{k}\sin k(t-\tau)\left(1-\sum\limits_{j = 1}^{\infty}\left(x_{j}^{1}(\tau, \varepsilon)+c_{1}^{j}\cos j\tau + c_{2}^{j}\sin j\tau\right)^2 \right)\times
\end{equation*}
\begin{equation*}
\times\left(y_{k}^{1}(\tau, \varepsilon) - c_{1}^{k}\sin k\tau + c_{2}^{k}\cos k\tau\right),
\end{equation*}
\begin{equation*}
\overline{H}_{2k}(t, \tau, \varepsilon) = \frac{1}{k}\cos k(t-\tau)\left(1-\sum\limits_{j = 1}^{\infty}\left(x_{j}^{1}(\tau, \varepsilon)+c_{1}^{j}\cos j\tau + c_{2}^{j}\sin j\tau\right)^2 \right)\times
\end{equation*}
\begin{equation*}
\times\left(y_{k}^{1}(\tau, \varepsilon) - c_{1}^{k}\sin k\tau + c_{2}^{k}\cos k\tau\right).
\end{equation*}

Substituting the representation (\ref{eqva:4}) in the condition for the solvability (\ref{eqva:3}), we obtain the operator equation (\ref{Pok-hyperbolic-eqcut:6}) with the operator $B_{0}$, which, in this case, has the form:
\begin{equation*}
B_{0} = -\int\limits_{0}^{2\pi}U^{-1}(\tau)\left( A_1(\tau)U(\tau) + A_2(\tau)U(\tau)\right)d\tau
\end{equation*}
and the element $b$, which is defined as follows:
\begin{equation*}
b = \int\limits_{0}^{2\pi}U^{-1}(\tau)\left(A_1(\tau)\overline{z}(\tau, \varepsilon) + A_{2}(\tau)\overline{z}(\tau, \varepsilon) + R(z(\tau, \varepsilon)) \right)d\tau
\end{equation*}
or
\begin{align}
b = {}& \varepsilon \int\limits_{0}^{2\pi}U^{-1}(\tau)\left(A_1(\tau) + A_2(\tau)\right)\int\limits_{0}^{\tau}\overline{H}(\tau, s, \varepsilon)dsd\tau + \int\limits_{0}^{2\pi} U^{-1}(\tau)R(z(\tau, \varepsilon))d\tau \notag\\
 = {}& \varepsilon \int\limits_{0}^{2\pi}U^{-1}(\tau)\left(A_1(\tau) + A_2(\tau)\right)\int\limits_{0}^{\tau}U(\tau -s)H\left(z(s,\varepsilon) + z_{0}\left(s, c^{0}\right)\right)dsd\tau  \notag \\
{}& + \int\limits_{0}^{2\pi} U^{-1}(\tau)R(z(\tau, \varepsilon))d\tau = \varepsilon b_1 + b_2. \label{ura:3}
\end{align}
The components of the vectors
\begin{align*}
b_1 = {}& {\rm col}\begin{pmatrix} b_{11}, & b_{12}, & \cdots, & b_{1(2k-1)}, & b_{1(2k)}, & \cdots \end{pmatrix}, \\
b_2 = {}& {\rm col} \begin{pmatrix} b_{21}, & b_{22}, & \cdots &, b_{2(2k-1)}, & b_{2(2k)}, & \cdots \end{pmatrix}
\end{align*}
have the form:
\begin{equation*}
b_{1(2k-1)} = -\sum_{p = 1}^{N}\frac{2}{k}\int\limits_{0}^{2\pi} \sin k\tau x_{p}^{0}\left(\tau, c^{0}\right)y_{k}^{0}\left(\tau,c^{0}\right)\times
\end{equation*}
\begin{equation*}
\times \int\limits_{0}^{\tau}\sin p(\tau - s)\left(1 - \sum_{j = 1}^{\infty}\left(x_{j}^{0}\left(s, c^{0}\right)\right)^{2}\right)\left(y_{p}^{1}(s,\varepsilon) - y_{p}^{0}\left(s, c^{0}\right)\right)dsd\tau -
\end{equation*}
\begin{equation*}
- \frac{1}{k}\int\limits_{0}^{2\pi}\sin k\tau \left(1 - \sum_{j = 1}^{N}\left(x_{j}^{0}\left(\tau, c^{0}\right)\right)^{2} \right)\times
\end{equation*}
\begin{equation*}
\times \int\limits_{0}^{\tau}\cos k(\tau - s)\left(1 - \sum_{j = 1}^{\infty}\left(x_{j}^{0}\left(s, c^{0}\right)\right)^{2} \right)\left(y_{k}^{1}(s, \varepsilon) - y_{k}^{0}\left(s, c^{0}\right)\right)dsd\tau,
\end{equation*}
\begin{equation*}
b_{1(2k)} = \sum_{p = 1}^{N}\frac{2}{k}\int\limits_{0}^{2\pi} \cos k\tau x_{p}^{0}\left(\tau, c^{0}\right)y_{k}^{0}\left(\tau,c^{0}\right)\times
\end{equation*}
\begin{equation*}
\times \int\limits_{0}^{\tau}\sin p(\tau - s)\left(1 - \sum_{j = 1}^{\infty}\left(x_{j}^{0}\left(s, c^{0}\right)\right)^{2}\right)\left(y_{p}^{1}(s,\varepsilon) - y_{p}^{0}\left(s, c^{0}\right)\right)dsd\tau -
\end{equation*}
\begin{equation*}
+ \frac{1}{k}\int\limits_{0}^{2\pi}\cos k\tau \left(1 - \sum_{j = 1}^{N}\left(x_{j}^{0}\left(\tau, c^{0}\right)\right)^{2} \right)\times
\end{equation*}
\begin{equation*}
\times \int\limits_{0}^{\tau}\cos k(\tau - s)\left(1 - \sum_{j = 1}^{\infty}\left(x_{j}^{0}\left(s, c^{0}\right)\right)^{2} \right)\left(y_{k}^{1}(s, \varepsilon) - y_{k}^{0}\left(s, c^{0}\right)\right)dsd\tau,
\end{equation*}
$k = \overline{1, N}$ and $b_{1(s)} = 0$, $s \geq 2N$.

Similarly
\begin{align*}
b_{2(2k-1)} = {}& \frac{1}{k}\int\limits_{0}^{2\pi}\sin k\tau\left(\sum\limits_{j = 1}^{\infty}\left(x_{j}^{1}(\tau, \varepsilon)\right)^2\left(y_{k}^{0}\left(\tau, c^{0}\right) + y_{k}^{1}(\tau, \varepsilon)\right) \right. \\
  {}& \left. + \sum\limits_{j = 1}^{N}x_{j}^{0}\left(\tau, c^{0}\right)x_{j}^{1}(\tau, \varepsilon)y_{k}^{1}(\tau, \varepsilon)\right)d\tau, \\
b_{2(2k)} = {}& - \frac{1}{k}\int\limits_{0}^{2\pi}\cos k\tau\left(\sum\limits_{j = 1}^{\infty}\left(x_{j}^{1}(\tau, \varepsilon)\right)^2\left(y_{k}^{0}\left(\tau, c^{0}\right) + y_{k}^{1}(\tau, \varepsilon)\right) \right. \\
  {}& \left. + \sum\limits_{j = 1}^{N}x_{j}^{0}\left(\tau, c^{0}\right)x_{j}^{1}(\tau, \varepsilon)y_{k}^{1}(\tau, \varepsilon)\right)d\tau.
\end{align*}

Note, in this case, the operator $B_{0}$ is a countable-dimensional matrix and can be found using the operator $F$ (\ref{porodg:1}) as follows:
\begin{equation} \label{eq:34}
B_{0} = \frac{dF\left(c^{0}\right)}{dc} = \begin{pmatrix} \frac{\partial F_{1}^{1}}{\partial c_1^{1}} & \frac{\partial F_{1}^{1}}{\partial c_2^{1}} & \cdots & \frac{\partial F_{1}^{1}}{\partial c_{1}^{k}} & \frac{\partial F_{1}^{1}}{\partial c_{2}^{k}} & \cdots  \\
\frac{\partial F_{1}^{2}}{\partial c_1^{1}} & \frac{\partial F_{1}^{2}}{\partial c_2^{1}} & \cdots & \frac{\partial F_{1}^{2}}{\partial c_{1}^{k}} & \frac{\partial F_{1}^{2}}{\partial c_{2}^{k}} & \cdots \\
\cdots & \cdots & \cdots & \cdots & \cdots & \cdots \\
\frac{\partial F_{k}^{1}}{\partial c_1^{1}} & \frac{\partial F_{k}^{1}}{\partial c_2^{1}} & \cdots & \frac{\partial F_{k}^{1}}{\partial c_{1}^{k}} & \frac{\partial F_{k}^{1}}{\partial c_{2}^{k}} & \cdots \\
\frac{\partial F_{k}^{2}}{\partial c_1^{1}} & \frac{\partial F_{k}^{2}}{\partial c_2^{1}} & \cdots & \frac{\partial F_{k}^{2}}{\partial c_{1}^{k}} & \frac{\partial F_{k}^{2}}{\partial c_{2}^{k}} & \cdots \\
\cdots & \cdots & \cdots & \cdots & \cdots & \cdots  \end{pmatrix} = \left(\frac{dF}{dc}\right)_{k,j=1}^{\infty}.
\end{equation}
Using (\ref{porodg:01}), (\ref{porodg:001}), we get:
\begin{align*}
{}& \frac{\partial F_{1}^{k}}{\partial c_{1}^{k}} = -\frac{\pi}{2k}\left(c_{1}^{k0}\right)^{2}, \quad
\frac{\partial F_{2}^{k}}{\partial c_{2}^{k}}= -\frac{\pi}{2k} \left(c_{2}^{k0}\right)^{2}, \\
{}& \frac{\partial F_{1}^{k}}{\partial c_{2}^{k}} = \frac{\partial F_{2}^{k}}{\partial c_{1}^{k}} = -\frac{\pi}{2k}c_{1}^{k0}c_{2}^{k0}, \quad j=k,
\end{align*}
and
\begin{align*}
{}& \frac{\partial F_{1}^{k}}{\partial c_{1}^{j}} = -\frac{\pi}{k}c_{1}^{k0}c_{1}^{j0}, \quad \frac{\partial F_{1}^{k}}{\partial c_{2}^{j}} = -\frac{\pi}{k}c_{1}^{k0}c_{2}^{j0}, \\
{}& \frac{\partial F_{2}^{k}}{\partial c_{1}^{j}} = -\frac{\pi}{k}c_{2}^{k0}c_{1}^{j0}, \quad \frac{\partial F_{2}^{k}}{\partial c_{2}^{j}} = -\frac{\pi}{k}c_{2}^{k0}c_{2}^{j0}, \quad j\neq k.
\end{align*}
Therefore, the matrix $B_{0}$ has the following form:
\begin{equation*}
B_{0} = -\frac{\pi}{2} \begin{pmatrix}
B^{0}_{11} & B^{0}_{12} & \cdots & B^{0}_{1N} & O & \cdots  \\
B^{0}_{21} & B^{0}_{22} & \cdots & B^{0}_{2N} & O & \cdots  \\
\vdots & \vdots  & \ddots & \vdots & \vdots  & \vdots  \\
B^{0}_{k1} & B^{0}_{k2} & \cdots & B^{0}_{NN} & O & \cdots  \\
O & O & \cdots & O & O & \cdots  \\
\cdots  & \cdots  & \cdots & \cdots  & \cdots & \cdots
\end{pmatrix},
\end{equation*}
\begin{equation*}
  B^{0}_{kj}= \left\{
         \begin{array}{ll}
           \displaystyle\frac{1}{k}
\begin{pmatrix}
\left(c_{1}^{k0}\right)^{2} & c_{1}^{k0}c_{2}^{k0} \\
c_{1}^{k0}c_{2}^{k0} & \left(c_{2}^{k0}\right)^{2} \
\end{pmatrix}, & \hbox{$k=j$;} \\
           \displaystyle\frac{2}{k}
\begin{pmatrix}
c_{1}^{k0}c_{1}^{j0} & c_{1}^{k0}c_{2}^{j0} \\
c_{2}^{k0}c_{1}^{j0} & c_{2}^{k0}c_{2}^{j0} \
\end{pmatrix}, & \hbox{$k\neq j$.}
         \end{array}
       \right.
\end{equation*}

Let us check the solvability condition (\ref{Pok-hyperbolic-eqcut:7}) of the boundary-value problem (\ref{Pok-hyperbolic-eqde:11}), (\ref{Pok-hyperbolic-eqde:12}). At first we find a pseudo-inverse matrix $B_{0}^{+}$. Note, that the matrix $B_{0}$ can be represented as the product of three infinite-dimensional matrices, i.e.
\begin{equation*}
B_{0} = V_{1} W V_{2},
\end{equation*}
where
\begin{equation*}
W = -\frac{\pi}{2} \begin{pmatrix}
W_{1} & W_{2} & W_{2} & \cdots & W_{2} & O & \cdots  \\
W_{2} & W_{1} & W_{2} & \cdots & W_{2} & O & \cdots  \\
W_{2} & W_{2} & W_{1} & \cdots & W_{2} & O & \cdots  \\
\vdots & \vdots & \vdots  & \ddots & \vdots & \vdots  & \vdots  \\
W_{2} & W_{2}  & W_{2} & \cdots & W_{1} & O  & \cdots  \\
O & O  & O & \cdots & O & O  & \cdots  \\
\cdots  & \cdots  & \cdots  & \cdots & \cdots & \cdots  & \cdots
\end{pmatrix},
\end{equation*}
\begin{equation*}
V_{1}=  {\rm diag} \begin{pmatrix}
V^{1}_{1} & V^{1}_{2} & \cdots & V^{1}_{N} & O & \cdots
\end{pmatrix},
\end{equation*}
\begin{equation*}
V_{2}=  {\rm diag} \begin{pmatrix}
V^{2}_{1} & V^{2}_{2} & \cdots & V^{2}_{N} & O & \cdots
\end{pmatrix},
\end{equation*}
\begin{equation*}
V^{1}_{k} =\frac{1}{k}
\begin{pmatrix}
c_{1}^{k0} & 0 \\
0 & c_{2}^{k0}
\end{pmatrix}, \quad
V^{2}_{k} =
\begin{pmatrix}
c_{1}^{k0} & 0 \\
0 & c_{2}^{k0}
\end{pmatrix},
\quad
W_{k} = k
\begin{pmatrix}
1 & 1 \\
1 & 1
\end{pmatrix}.
\end{equation*}

In our case, the diagonal matrices $V_{1}$, $V_{2}$ and the symmetric matrix $W$ have non-zero blocks of size $2N\times 2N$. According to \cite{Gouveia-Puystjens-1991}, the matrix $B_{0}^{+}$ can be found by one of the following formulas:
\begin{align}
B_{0}^{+} = {}& (WV_{2})^{+}W(V_{1}W)^{+}, \label{eq-psevdo:3} \\
B_{0}^{+} = {}& (WV_{2})^{*}[WV_{2}(WV_{2})^{*}+I-WW^{+}]^{-1} \notag\\
{}& W[(V_{1}V_{2})^{*}V_{1}W+I-W^{+}W]^{-1}(V_{1}W)^{*}.\label{eq-psevdo:4}
\end{align}
The matrix $W$ has a special structure, so the pseudo-inverted matrix $W^{+}$ from the formula (\ref{eq-psevdo:4}) can be found explicitly
\begin{equation}\label{eq-psevdo:1}
W^{+} = -\frac{1}{2(2N-1)\pi} \begin{pmatrix}
W_{3-2N} & W_{2} & \cdots & W_{2} & O & \cdots  \\
W_{2} & W_{3-2N} & \cdots & W_{2} & O & \cdots  \\
W_{2} & W_{2} & \cdots & W_{2} & O & \cdots  \\
\vdots & \vdots  & \ddots & \vdots & \vdots  & \vdots  \\
W_{2} & W_{2} & \cdots & W_{3-2N} & O  & \cdots  \\
O & O & \cdots & O & O  & \cdots  \\
\cdots  & \cdots  & \cdots & \cdots & \cdots  & \cdots
\end{pmatrix}.
\end{equation}

\begin{rmk}\label{Pok-hyperbolic-rem-2}
The matrix $B^{-}_{0}= V^{-1}_{2} W V^{-1}_{1}$ is a generalized inverse matrix to the matrix $B_{0}$.
\end{rmk}

The condition (\ref{Pok-hyperbolic-eqcut:7}), in the general case, is not true. Therefore, the expression (\ref{eqva:4}) defines the pseudosolutions of the periodic boundary-value problem for the abstract Van der Pol equation. Thus, to find a pseudo-solution of this problem, we can apply the iterative process of the Theorem \ref{Pok-hyperbolic-th-3} and we obtain the following corollary.

\begin{crl}\label{Pok-hyperbolic-cor-1} If the pairs of constants $(c_{1}^{k0}, c_{2}^{k0})$, $k = \overline{1, N}$  satisfy all the conditions of the Theorem \ref{Pok-hyperbolic-th-5} and the condition (\ref{Pok-hyperbolic-eq:7}), then the periodic boundary-value problem for the abstract Van der Pol equation (\ref{Pok-hyperbolic-eqde:11}), (\ref{Pok-hyperbolic-eqde:12}) has a pseudosolution  and can be found by the following iterative process:
\begin{align*}
z^{m + 1}(t, \varepsilon) = {}& U(t)c^{m}+ \overline{z}^{m+1}(t, \varepsilon), \\
c^{m } = {}& B_{0}^{+} b^{m}  = B_{0}^{+}\left(\varepsilon b_1^m + b_{2}^{m}\right), \\
\overline{z}^{m + 1}(t, \varepsilon) = {}& \varepsilon\int\limits_{0}^{t}U(t - \tau)H\left(z^{m}(\tau, \varepsilon) + z_{0}\left(\tau, c^{0}\right)\right)d\tau,
\end{align*}
where the matrix $B_{0}^{+}$ is determined by the relations (\ref{eq-psevdo:3}), (\ref{eq-psevdo:4}), the vectors $b$, $b_1$, $b_2$ have the form (\ref{ura:3}) and the vector function $\overline{z}(t, \varepsilon)$ is determined by the relation (\ref{ura:2}).
\end{crl}

\section{Conclusions}

In this paper we develop constructive methods of investigation of boundary-value problems for the hyperbolic equation in the Hilbert and Banach spaces. For a linear problem we find a necessary and sufficient conditions of strong generalized solvability which are analogeus to the Fredholm alternative. For a weakly nonlinear problem we find conditions of branching of solutions. As an application we construct theory of generalized pseudoinvertibility for van der Pol equation. Proposed in the article approaches gives possibility  to study boundary-value problems in the infinite-dimensional spaces from a single point of view. 





\bibliography{Boichuk_Pokutnyi_Feruk_bibfile}

\end{document}